\documentclass[12pt]{amsart}
\usepackage{amssymb,latexsym}
\usepackage{enumerate}
\usepackage{amscd}
\usepackage{verbatim}

\theoremstyle{plain}
\newtheorem{theorem}{Theorem}

\newtheorem{lemma}[theorem]{Lemma}

\theoremstyle{definition}

\theoremstyle{remark}

\def\pair#1#2{\langle #1, #2\rangle}

\def\Mof#1{\Cal M(A)}

\newcommand{\Con}{\operatorname{Con}}

\newcommand{\nat}{\operatorname{nat}}

\def\congruence{on\-gru\-ence\discretionary{-}{}{-}}
\def\conM/{c\congruence mod\-u\-lar}
\def\ConM/{C\congruence mod\-u\-lar}
\def\conD/{c\congruence dis\-trib\-u\-tive}
\def\ConD/{C\congruence dis\-trib\-u\-tive}
\def\conP/{c\congruence per\-mut\-a\-ble}
\def\ConP/{C\congruence per\-mut\-a\-ble}
\def\conMity/{\conM/\-i\-ty}
\def\ConMity/{\ConM/\-i\-ty}
\def\conDity/{c\congruence dis\-trib\-u\-tiv\-i\-ty}
\def\ConDity/{C\congruence dis\-trib\-u\-tiv\-i\-ty}
\def\conPity/{c\congruence per\-mut\-a\-bil\-i\-ty}
\def\ConPity/{C\congruence per\-mut\-a\-bil\-i\-ty}
\def\usprv/{un\-der\-ly\-ing-set-pre\-ser\-ving}

\def\ie/{{i.e.}}
\def\Ie/{{I.e.}}
\def\eg/{{e.g.}}
\def\Eg/{{E.g.}}
\def\etc/{{etc.}}

\newdimen\mysubdimen
\newbox\mysubbox

\def\subwhat#1#2#3{{
\setbox\mysubbox=\hbox{#3}
\mysubdimen=\wd\mysubbox
\setbox\mysubbox=\hbox{$#1#2$}
\ifnum\mysubdimen>\wd\mysubbox
\vtop{
\hbox to\mysubdimen{\hfil\box\mysubbox\hfil}
\nointerlineskip
\hbox{#3}}
\else
\mysubdimen=\wd\mysubbox
\vtop{
\box\mysubbox
\nointerlineskip
\hbox to\mysubdimen{#3}}
\fi
}}
\begin{document}

% One author
\title{J\'onsson's Theorem in Non-modular Varieties}
\author{William H. Rowan}
\address{PO Box 20161 \\
         Oakland, California 94620--0161}
\email{whrowan@member.ams.org}
%\thanks{thanks}
% End one author

%\def\cite#1{[#1]}

\keywords{centralizer, commutator, J\'onsson's Theorem}
\subjclass{Primary: 08B99; Secondary: 03C05, 03C20,
08B10, 08B26}
\date{\today}

\begin{abstract}
A version of J\'onsson's Theorem, as previously
generalized, holds in non-modular varieties.
\end{abstract}
\maketitle

\section*{Introduction}

In this paper, we state and prove a further
generalization of J\'onsson's Theorem which holds in
non-modular varieties. Our starting point is
\cite[Theorem 10.1]{F-McK}, a generalization of
J\'onsson's Theorem due to Hagemann, Hermann, Freese,
McKenzie, and Hrushovskii:

\begin{theorem}
Let $\mathcal K$ be a class of algebras, and suppose that
$\mathbf{HSP}(\mathcal K)$ is congruence-modular. If
$B\in\mathbf{HSP}(\mathcal K)$ is subdirectly irreducible
and $\alpha$ is the centralizer of the monolith of $B$,
then $B/\alpha\in\mathbf{HSP}_u(\mathcal K)$.
\end{theorem}

Here $\mathbf H$, $\mathbf S$, and $\mathbf P$ are the
standard operators on classes of algebras that close them
under homomorphic images, subalgebras, and products
respectively, and $\mathbf P_u$ is the closure under
ultraproducts. Also, recall that the \emph{monolith}
$\mu$ of a subdirectly irreducible algebra $B$ is the
minimal nontrivial congruence of $B$, and the
\emph{centralizer} of $\mu$ is the largest congruence
$\alpha$ such that $[\mu,\alpha]=\bot$, which exists
because of the additivity of the commutator operation.

In \S1, we introduce a notion of centrality
which we call \emph{lax centrality}, which
reduces, as shown in \S2, to the usual,
commutator-theoretic notion in modular varieties. This
definition was motivated by the requirements of the proof
of the above theorem, and by its somewhat
tenuous connection, discussed in \S3, with the concept of
the
\emph{free intersection},
which generalizes the modular commutator to non-modular
varieties \cite{D-G} \cite{D-K}. In
\S4, we give the further
generalization of Theorem~1 that is the goal of the paper.
We make some closing remarks in \S5.

\section*{Preliminaries}

\subsection{Category Theory} We follow
\cite{MacL} in terminology and notation.

\subsection{Lattices}
We use the symbols $\bot$ and $\top$ to denote the least
and greatest elements of a lattice, assuming they exist.

\subsection{Universal Algebra} We assume the basic
definitions of Universal Algebra, as found, for example,
in
\cite{B-S}. However, in the definition of an algebra, we
prefer to allow an algebra to be empty.

In \S2, we will assume an acquaintance with Commutator
Theory for congruence-modular varieties, as developed in
\cite{F-McK}. Note that the notation used in \cite{F-McK}
(and some other works on Commutator Theory) is different
from ours. In particular, the least and greatest elements
of a congruence lattice are denoted by
$0$ and $1$ rather than $\bot$ and $\top$.

If $A$, $B$ are algebras, $f:A\to B$ is a homomorphism,
and $\alpha\in\Con A$, we denote by $\vec f\alpha$ the
congruence of $B$ generated by $f(\alpha)$, the set of
pairs $\pair{f(a)}{f(a')}$ such that $a\mathrel\alpha
a'$.

\begin{theorem} Let $A$, $B$ be algebras and $f:A\to B$ a
homomorphism. Let $\alpha\in\Con A$. We have
\begin{enumerate}
\item If $f$ is onto, then $\vec f\alpha=f(\alpha\vee\ker
f)$;
\item if $\beta\in\Con B $, then $\vec f\alpha\leq\beta$
iff $\alpha\leq f^{-1}(\beta)$; and
\item if $C$ is another algebra and $g:B\to C$, then
$\vec{(gf)}\alpha=\vec g(\vec f\alpha)$.
\end{enumerate}
\end{theorem}

\section{Lax Centrality}

Let $\mathbf V$ be a variety of algebras of some type,
and $B\in\mathbf V$. Let $\mu\in\Con B$. We say that
$\alpha\in\Con B$ \emph{laxly centralizes} $\mu$ (with
respect to the variety $\mathbf V$) if there is an
algebra $C\in\mathbf V$, an onto homomorphism $\pi:C\to
B$, and $\beta$, $\gamma\in\Con C$ such that
$\vec\pi\beta\geq\mu$,
$\vec\pi\gamma\geq\alpha$, and
$\beta\wedge\gamma=\bot$.

Note that if $C$, $\pi$, $\beta$, and $\gamma$ are given
such that $\beta\wedge\gamma\leq\ker\pi$ and all of the
other conditions of the definition are satisfied, then
$\alpha$ laxly centralizes $\mu$, because we can replace
$C$ by $C/(\beta\wedge\gamma)$.

Also, note that if $\alpha$ laxly centralizes $\mu$,
$\alpha'\leq\alpha$, and $\mu'\leq\mu$, then $\alpha'$
laxly centralizes $\mu'$.

For some notions of centrality, we know that, given
$\mu\in\Con B$, there exists some maximal $\alpha$ such
that $\alpha$ centralizes $\mu$. However, we do not know
this for lax centrality.

\section{Lax Centrality and the Modular Commutator}

The following theorem shows that if $\mathbf V$ is
congruence-modular, the notion of lax centrality is
identical to the notion of centrality given by
commutator theory.

\begin{theorem}\label{T:Modular} Let $\mathbf V$ be a
congruence-modular variety of algebras, and let
$B\in\mathbf V$ and
$\mu$, $\alpha\in\Con B$. Then $\alpha$ laxly
centralizes
$\mu$ iff
$[\mu,\alpha]=\bot$.
\end{theorem}

\begin{proof} If $\alpha$ laxly centralizes $\mu$, let
$C$, $\pi$, $\beta$, and $\gamma$ be given as in the
definition of lax centrality.  Let $\ker\pi$ be
denoted by $\theta$.
Since
$\beta\wedge\gamma=\bot$, we have $[\beta,\gamma]=\bot$,
which implies that
$[\vec\pi\beta,\vec\pi\gamma]=
[\mu,\alpha]=\bot$.

If $[\mu,\alpha]=\bot$, we let $C=B(\mu)$ be the
subalgebra of $B\times B$ of pairs related by $\mu$, and
form the pushout square
\[
\begin{CD} C @<<\Delta< B \\
@V\rho VV @VV\nat\alpha V \\
P @<<< B/\alpha
\end{CD}
\]
where $\Delta$ is defined by $\Delta(b)=\pair bb$.
Let
$\pi:C\to B$ be given by $\pi:\pair bc\mapsto c$,
$\beta=\ker\pi_\mu$ where
$\pi_\mu:C\to B$ is given by $\pi_\mu:\pair bc\mapsto b$,
and $\gamma=\ker\rho$ where $\rho$ is the homomorphism
shown in the pushout square.  Note that $\rho$ is onto,
because it is a pushout of the onto homomorphism
$\nat\alpha$.  In the notation of \cite[Chapter
IV]{F-McK}, we have $\gamma=\Delta_{\mu,\alpha}$.

Let $\theta = \ker\pi$.

We have $\vec\pi\beta=\mu$, and we have
$\vec\pi\rho=\alpha$ because that is the
pushout of $\rho$ along $\pi$, and $\pi\Delta=1_B$.

To show that $\beta\wedge\gamma=\bot$, we must use the
fact that
$[\mu,\alpha]=\bot$. Indeed, by \cite[Theorem
4.9(iv)]{F-McK}, we have
\begin{align*}
\beta\wedge\gamma
&= \ker\pi_\mu\wedge\Delta_{\mu,\alpha} \\
&= \{\,\pair{\pair xy}{\pair zw}\mid x=z\hbox{ and }
y\mathrel{[\mu,\alpha]}w\,\}
\\
&= \bot.
\end{align*}

\end{proof}

\section{Lax Centrality and the Free Intersection}

Let $B$ be an algebra in a variety $\mathbf V$, and
$\alpha$, $\beta\in\Con B $. Let $F$ be the relatively
free algebra in $\mathbf V$ on the set of generators
$\{\,x_b,y_b\,\}_{b\in B}$, i.e., on the disjoint union
of two copies of the underlying set of $B$. The
\emph{free intersection} of $\alpha$ and $\beta$
(with respect to $\mathbf V$) is defined \cite{D-G}
\cite{D-K} as
$\vec\zeta(\bar\alpha\wedge\bar\beta)$, where $\bar\alpha$
is the congruence on $F$ generated by the relation
$\{\,\pair{x_b}{x_{b'}}\mid b\mathrel\alpha b'\,\}$ on the
generators of $F$, $\bar\beta$ is the congruence
generated by the relation $\{\,\pair{y_b}{y_{b'}}\mid
b\mathrel\beta b'\,\}$, and $\zeta$ is the onto
homomorphism defined on generators by $x_b\mapsto b$,
$y_b\mapsto b$. Note that $\vec\zeta\bar\alpha=\alpha$ and
$\vec\zeta\bar\beta=\beta$.

We will denote the free intersection of $\alpha$ and
$\beta$ by $[\alpha,\beta]$.  This use of commutator
notation is justified by the fact (\cite[Theorem
5.14]{D-K} or \cite[Theorem 2.4]{D-G}) that in modular
varieties, the free intersection of two congruences is
their commutator. Clearly, we have
$[\alpha,\beta]=[\beta,\alpha]$, and the free
intersection is monotonic in its arguments.

Let $B$ belong to a variety $\mathbf V$ of
algebras of some type, and let
$\mu$,
$\alpha\in\Con B$. Clearly, if
$[\mu,\alpha]=\bot$, then $\alpha$ laxly centralizes
$\mu$.

If the free intersection $[\mu,\alpha]$ is the
minimum congruence of
$B$ of the form $\vec\pi(\bar\mu\wedge\bar\alpha)$, where
$\pi$ is onto, $\vec\pi\bar\mu=\mu$, and
$\vec\pi\bar\alpha=\alpha$, a proof of which is not
known to us, that would imply the converse, i.e., that if
$\alpha$ laxly centralizes $\mu$, then
$[\mu,\alpha]=\bot$.

\section{Lax Centrality and J\'onsson's Theorem}

We begin with a lemma:

\begin{lemma} Let $\mathcal K$ be a class of algebras of
some type, and let $B\in\mathbf{HSP}(\mathcal K)$ and
$\mu$, $\alpha\in\Con B$. If $\alpha$ laxly centralizes
$\mu$, then there exist $C$, $\pi$, $\beta$, and $\gamma$,
as in the definition of lax centrality, such that in
addition we have $C/\beta$,
$C/\gamma\in\mathbf{SP}(\mathcal K)$.
\end{lemma}

\begin{proof}
Let $C_0$, $\pi_0$, $\beta_0$, and $\gamma_0$ witness
the fact that $\alpha$ laxly centralizes $\mu$. Since
$\beta_0\wedge\gamma_0=\bot$, there is a natural one-one
homomorphism $\iota_0:C_0\to C_0/\beta_0\times
C_0/\gamma_0$.

It is
assumed in the definition of lax centrality that
$C_0\in\mathbf{HSP}(\mathcal K)$, and the same goes for
$C_0/\beta_0$ and $C_0/\gamma_0$. Thus, there are
$C_\beta$, $C_\gamma\in\mathbf{SP}(\mathcal K)$ and
onto homomorphisms $\pi_\beta:C_\beta\to C_0/\beta_0$,
$\pi_\gamma:C_\gamma\to C_0/\gamma_0$.

Let $C=\{\,\langle c,b,g\rangle\mid c\in C_0,\,
b\in C_\beta,\,g\in C_\gamma,\,c/\beta_0=\pi_\beta(b),
\text{ and }c/\gamma_0=\pi_\gamma(g)\,\}$. Then
the diagram
\[
\begin{CD} C @>\iota_1>> C_\beta\times C_\gamma \\
@V\pi_1 VV @VV\pi_\beta\times\pi_\gamma V \\
C_0 @>>\iota_0 > C_0/\beta_0\times C_0/\gamma_0
\end{CD}
\]
is a pullback square, where $\pi_1\langle c,b,g\rangle =
c$ and $\iota_1\langle c,b,g\rangle = \langle b,g\rangle$.
Clearly $\pi_1$ is onto, and since $\iota_0$ is one-one,
$\iota_1$ is one-one as well.

Let $\beta$ be the kernel of the
homomorphism $\nu_\beta:C\to C_\beta$ given by
$\nu_\beta\langle c,b,g\rangle = b$, and $\gamma$
the kernel of the homomorphism $\nu_\gamma:C\to C_\gamma$
given by
$\nu_\gamma\langle c,b,g\rangle=g$. Then
$\beta\wedge\gamma=\bot$, because
$\iota_1$ is one-one.  Also, $\nu_\beta$ and $\nu_\gamma$
are easily seen to be onto.  Thus, $C/\beta$,
$C/\gamma\in\mathbf{SP}(\mathcal K)$.

Let $\pi=\pi_0\pi_1$. To show that $C$,
$\pi$,
$\beta$, and
$\gamma$ witness the lax centrality of $\mu$ by
$\alpha$, it suffices to prove that
$\vec\pi_1\beta=\beta_0$ and
$\vec\pi_1\gamma=\gamma_0$.

Suppose $c\mathrel{\beta_0}c'$. There is a $b\in C_\beta$
such that $\pi_\beta(b)=c/\beta_0=c'/\beta_0$.
Let $g$, $g'\in C_\gamma$ be
such that $\pi_\gamma(g)=c/\gamma_0$ and
$\pi_\gamma(g')=c'/\gamma_0$. Then $\langle
c,b,g\rangle\mathrel\beta\langle c',b,g'\rangle$. It
follows that $\vec\pi_1\beta\geq\beta_0$.

On the other hand,
$(\pi_1)^{-1}(\beta_0)=\{\,\langle\langle
c,b,g\rangle,\langle c',b',g'\rangle\rangle\in C^2\mid
c\mathrel{\beta_0}c'\,\}$. If $\langle
c,b,g\rangle\mathrel\beta\langle c',b',g'\rangle$, then
we have $b=b'$, which implies that $c\mathrel\beta_0 c'$.
Thus, we have $\beta\leq(\pi_1)^{-1}(\beta_0)$. It follows
that $\vec\pi_1\beta\leq\beta_0$.

That $\vec\pi_1\gamma=\gamma_0$ follows by
similar arguments.
\end{proof}

Now for the main theorem:

\begin{theorem} Let $\mathcal K$ be a class of algebras
of some type, and let
$B$ be a subdirectly irreducible algebra in
$\mathbf{HSP}(\mathcal K)$, with monolith $\mu$. Let
$\alpha\in\Con B$ be maximal for the property that
$\alpha$ laxly centralizes $\mu$. Then $B/\alpha\in
\mathbf{HSP}_u(\mathcal K)$.
\end{theorem}

\begin{proof} Let $C$, $\pi$, $\beta$, and $\gamma$ be
given as in the definition of lax centrality, and
such that $C/\beta$, $C/\gamma\in\mathbf{SP}(\mathcal K)$,
and denote
$\ker\pi$ by
$\theta$.

Let
$\{A_i\}_{i\in I}$, $\{A_i\}_{i\in \bar I}$ be tuples of
algebras in $\mathcal K$, with $I$ and $\bar I$ disjoint,
such that $C/\beta\hookrightarrow\Pi_IA_i$ and
$C/\gamma\hookrightarrow\Pi_{\bar I}A_i$.
Since $\beta\wedge\gamma=\bot$, we have embeddings
$C\hookrightarrow C/\beta\times
C/\gamma\hookrightarrow\Pi_{I\cup\bar I}A_i$.

Let $\mathcal F$ be a filter on $I\cup\bar I$ maximal with
respect to the property that $J\in\mathcal F$ implies
$\eta_J\leq\theta$, where $\eta_J$ is the kernel of the
natural map from $C$ to $\Pi_{j\in J}A_j$. We have
$I\notin\mathcal F$, because $\beta\not\leq\theta$.

Let $\mathcal U$ be an ultrafilter on $I\cup\bar I$
extending $\mathcal F$ and containing $\bar I$. Such an
ultrafilter exists, because $I\notin\mathcal F$. We
claim that if $J\in\mathcal U$, then
$\eta_J\leq\pi^{-1}(\alpha)$. For, if
$J\in\mathcal F$, then $\eta_J\leq\theta$. On the other
hand, if
$J\in\mathcal U-\mathcal F$, we have by the maximality of
$\mathcal F$ that neither $J$ nor its complement $\bar J$
can be adjoined to $\mathcal F$. This implies that there
is a
$K\in\mathcal F$ such that $\eta_{J\cap
K}\not\leq\theta$ and $\eta_{\bar J\cap K}\not\leq\theta$.
But $\eta_{J\cap K}\wedge\eta_{\bar J\cap
K}=\eta_K\leq\theta$.
By the remark following the definition of lax
centrality,  $\vec\pi\eta_{J\cap K}$ laxly
centralizes $\vec\pi\eta_{\bar J\cap K}$. However,
$J\cap K\subseteq J$, so
$\eta_J\leq\eta_{J\cap K}$. Also,
$\mu\leq\vec\pi\eta_{\bar J\cap K}$, because $\mu$ is a
monolith. Thus,
$\vec\pi\eta_J$ laxly centralizes
$\mu$.
This holds for every $J\in\mathcal U$, so it
holds for $J\cap\bar I$ for any particular $J$.
We have $\vec\pi\eta_{J\cap\bar
I}\geq\vec\pi\eta_{\bar I}=\vec\pi\gamma\geq\alpha$ and
laxly centralizing $\mu$, but $\alpha$ is maximal for that
property. It follows that
for any
$J\in\mathcal U$,
$\vec\pi\eta_J\leq\vec\pi\eta_{J\cap\bar I}=\alpha$, which
verifies the claim.

Let $\eta_{\mathcal U}$ be the restriction to $C$ of the
ultrafilter congruence on $\Pi_{I\cup\bar I}A_i$. Since
$\eta_J\leq\pi^{-1}(\alpha)$ for all $J\in\mathcal U$, we
have
$\eta_{\mathcal U}\leq\pi^{-1}(\alpha)$. Thus,
$B/\alpha\in\mathbf{HSP}_u(\mathcal K)$.
\end{proof}

\section{Remarks}

In the main theorem, because we do not know whether the
set of $\alpha$ laxly centralizing $\mu$ admits maximal
elements, we had to assume that $\alpha$ is a maximal
element of that set. Note that maximality of $\alpha$
is essential, even in the proof of the
congruence-modular version of J\'onsson's Theorem. It
is not the fact that $\alpha$ centralizes $\mu$ that
implies that $B/\alpha\in\mathbf{HSP}_u(\mathcal K)$,
but the fact that $\alpha$ is maximal for that
property. The case where there is no such maximal
$\alpha$, assuming this can indeed happen, is a likely
point where further work may clarify the situation and
perhaps allow a further elaboration of the theorem.

Lax centrality for the variety $\mathbf V$ is a $\mathbf
V$-tuple of binary relations $\lambda_B$ on $\Con B$, for
$B\in\mathbf V$. If $\mathbf V$ is congruence-modular,
then $\alpha\mathrel{\lambda_B}\mu$ iff
$[\mu,\alpha]=\bot$. We might ask, for general $\mathbf
V$, is there a suitable commutator, i.e., a $\mathbf
V$-tuple of binary operations $\kappa_B$ on $\Con B$ for
$B\in\mathbf V$, such that $\alpha\mathrel{\lambda_B}\mu$
iff $\mu\mathrel{\kappa_B}\alpha=\bot$? There certainly
is one, but it may not be very suitable: define
\[
\mu\mathrel{\kappa_B}\alpha = 
\begin{cases}
\bot,      &\text{if $\alpha$ laxly centralizes $\mu$;}\\
\top,      &\text{otherwise.}
\end{cases}
\]

In case $\mathbf V$ is congruence-modular, we have
$[\mu,\alpha]\leq[\vec\pi\beta,\vec\pi\gamma]=\vec\pi
[\beta,\gamma]\leq\vec\pi(\beta\wedge\gamma)$, whenever
$\pi:C\to B$ and $\beta$, $\gamma\in\Con C$ are such that
$\vec\pi\beta\geq\mu$ and $\vec\pi\gamma\geq\alpha$.
On the other hand, if we take $C=B(\mu)$, $\pi$, $\beta$,
and $\gamma$ as in the proof of Theorem~3, we have
\begin{align*}
\vec\pi(\beta\wedge\gamma)
&= \vec\pi\{\,\pair{\pair xy}{\pair zw}\mid x=z\hbox{ and
} y\mathrel{[\mu,\alpha]}w\,\}
\\
&= [\mu,\alpha].
\end{align*}
Thus, we can define $[\mu,\alpha]$ to be the meet over
all $\langle C,\pi,\beta,\gamma\rangle$ (where
$\vec\pi\geq\mu$ and $\vec\pi\gamma\geq\alpha$) of
$\vec\pi(\beta\wedge\gamma)$. Whether something similar
can be done for arbitrary $\mathbf V$ is an open
question; it may be that the meet can be $\bot$ even when
$\alpha$ does not laxly centralize $\mu$.

The referee rightly pointed out that we did not give an
effective method for computing an annihilator for $\mu$,
and that such a method would be helpful in effectively
applying Theorem~5.

More investigation will be required to fully
understand lax centrality, its place in Commutator
Theory, and its role in J\'onsson's Theorem.

\subsection*{Acknowledgement}
The author wishes to thank
Keith Kearnes for his careful reading of drafts of
this paper, for spotting an error in the original
proof of the main theorem, and for a number of useful
observations and references to the literature, and to
thank the referee for an insightful report.

\end{document}